\documentclass[12pt]{article}
\usepackage{amssymb}
\usepackage{amsfonts}
\usepackage{amsmath}
\usepackage[usenames]{color}
\usepackage{mathrsfs}
\usepackage{amsfonts}
\usepackage{amssymb,amsmath}
\usepackage{CJK}
\usepackage{cite}
\usepackage{cases}
\usepackage{amsthm}
\def\ssc{\scriptscriptstyle}
\def\cl{\centerline}

\def\al{\alpha}
\def\r{\gamma}
\def\b{\beta}
\def\vs{\vspace*}
\def\W{\mathscr{W}}
\def\L{\mathscr{W}_{L}(\Gamma){}}
\def\M{\mathscr{W}_{A}(\Gamma)}
\def\Z{\mathbb{Z}}

\def\C{\mathbb{C}}
\def\F{\mathbb{F}}
\def\QED{\hfill$\Box$}
\def\ni{\noindent}

\numberwithin{equation}{section}
\newtheorem{theo}{Theorem}[section]
\newtheorem{defi}[theo]{Definition}

\newtheorem{lemm}[theo]{Lemma}
\newtheorem{prop}[theo]{Proposition}

\begin{document}
\begin{center}
{\bf\large Lie bialgebras of generalized loop Virasoro algebras\,$^*$}
\footnote {$^*\,$Supported by NSF grant no.~10825101, 11001200, 11271284 and 11101269 of China, SMSTC grant
no.~12XD1405000, and the Fundamental Research Funds for the Central Universities.
\\\indent\ \ $^\dag\,$Corresponding author: X. Yue (xiaoqingyue@tongji.edu.cn).
}
\end{center}

\cl{Henan Wu, Song Wang, Xiaoqing Yue$^\dag$}

\cl{\small Department of Mathematics, Tongji University, Shanghai
200092, China}

\vs{8pt}

{\small
\parskip .005 truein
\baselineskip 3pt \lineskip 3pt

\noindent{{\bf Abstract:}
The first cohomology group of
a generalized loop Virasoro
algebra with coefficients in the tensor product of its adjoint module is shown to be trivial. The result is  applied to prove
that Lie bialgebra structures on generalized loop Virasoro
algebras are coboundary triangular. We then generalize the results to
generalized map Virasoro
algebras.\vs{5pt}

\ni{\bf Key words:}
Lie bialgebras, Yang-Baxter equation, generalized loop Virasoro algebras, generalized map Viarasoro algebras.}

\ni{\it Mathematics Subject Classification (2000):} 17B05, 17B37, 17B62, 17B68.}
\parskip .001 truein\baselineskip 6pt \lineskip 6pt
\section{Introduction}

As the universal central extension of the infinite dimensional Lie algebra (namely, the Witt algebra $W_1$)
of the linear differential operators $\{f(t)\frac d{dt}\,|\,f(t)\in\C[t,t^{-1}]\}$ of the Laurent polynomial algebra $\C[t,t^{-1}]$, the well-known Virasoro algebra plays a very fundamental role in conformal field theory, string theory, the theory of vertex operator
algebras  and the representation theory of
Kac-Moody algebras as well as extended affine Lie algebras (e.g., \cite{DL,Kac,AABGP}).
Various generalizations of the Virasoro algebra and other closely related algebras
have been objects of intensive studies in literature, e.g., \cite{GLZ,Sa,Su94,Su95,Su99,Su01,Su03,SWWY,WWY}. Among them, the one which naturally occurs in extended affine Lie algebras is the (generalized) loop Virasoro algebra (and more generally, the map Virasoro
algebras), whose structure and representation theories are studied  in \cite{GLZ,Sa,SWWY,WWY}.

In the present paper, we are interested in considering Lie
bialgebra structures on generalized loop Virasoro Lie algebras and generalized map Virasoro Lie algebras.
Lie bialgebras as well as their quantizations provide important tools in searching for solutions of quantum Yang-Baxter equations and in producing new quantum groups (e.g., \cite{D1,D2}).
Thus, there are a number of papers on Lie bialgebra structures  in literature. For instance,
Lie bialgebra
structures on the Witt, one-sided Witt and Virasoro algebras were shown in \cite{Ng,Taft} to be
coboundary triangular, which can be obtained from their nonabelian two dimensional Lie
subalgebras. Moreover, Lie bialgebra structures on the one-sided Witt algebra
were completely classified in \cite{Ng}. Furthermore, Lie bialgebra structures on
Lie algebras of generalized Witt type and Weyl type were considered
in \cite{SS,YS}.

We remark  that  the objects studied in this paper are $\Gamma$-graded Lie algebras (c.f.~\eqref{a1-0}), however, they are not finitely $\Gamma$-graded, and not
finitely generated as Lie algebras in general, in sharp contrast to
the Witt, one-sided Witt and Virasoro algebras.
Thus when we tackle  problems associated with such objects,
as stated in \cite{SWWY,WWY},
the classical techniques such as that in \cite{F} may not
be directly applied to the situation here. Instead, some new techniques or strategies must be employed as can be seen in
the proof of Lemma \ref{le3}. This is one of our motivations.
The main results presented in this paper are summarized in Theorems \ref{th1} and \ref{th2}, which state that all
Lie bialgebra structures under consideration are coboundary triangular.
As predicted, the results are not surprising. However, a Lie bialgebra having a
coboundary triangular structure does not mean that it is not interesting, as one can see from
\cite{SS1} that dualizing a coboundary triangular Lie bialgebra may produce new Lie algebras.
This will be persuaded in a sequel.

This paper proceeds as follows. In Section $2$, the preliminaries of Lie bialgebras are
recalled. In Section $3$, Lie bialgebra structures on generalized loop Virasoro algebras are
proven to be coboundary triangular. Then we generalize the results in Section $3$ to the case of map Virasoro algebras in Section $4$.

\section{Preliminaries}

Let us recall some concepts related to Lie bialgebras now.
Throughout this paper,
$\F$ denotes a field of characteristic zero and all the tensor products are taken over $\F$.
Let $L$ be a vector space over $\F$. Denote by $\tau$
the twist map of $L\otimes L$, namely,
\begin{equation}
\tau(x \otimes y)=y\otimes x, \ \ \forall\  x, y\in L.
\end{equation}
Denote by $\varepsilon$ the map which cyclically permutes the coordinates of $L\otimes L\otimes L$, i.e.,
\begin{equation}
\varepsilon(x_1\otimes x_2\otimes x_3)=x_2\otimes x_3\otimes x_1, \ \ \forall\  x_1, x_2, x_3\in L.
\end{equation}
Then the definition of a Lie algebra can be stated as follows.
\begin{defi}\rm
Let $L$ be a vector space over $\F$ and $\phi: L\otimes L \rightarrow L$ be a linear
map. The pair $(L, \phi)$ is called a Lie algebra if the following conditions are satisfied:
\begin{eqnarray*}
{\rm Ker}(1\otimes 1- \tau)\subset {\rm Ker}\phi \ \mbox{ and }\  \phi (1\otimes\phi)(1\otimes 1\otimes 1+\varepsilon+\varepsilon^2)=0,
\end{eqnarray*}
where $1$ denotes the identity map on $L$.
\end{defi}
The operator $\phi$ is usually called the bracket of $L$.
Dually, we get the definition of Lie coalgebras.
\begin{defi}\rm
Let $L$ be a vector space over $\F$ and $\Delta: L\rightarrow L\otimes L$ be a
linear map. The pair $(L, \Delta)$ is called a Lie coalgebra if
\begin{eqnarray*}
{\rm Im}\Delta \subset {\rm Im}(1\otimes 1-\tau)\ \mbox{ and }\  (1\otimes 1\otimes 1+\varepsilon+\varepsilon^2)(1\otimes \Delta)\Delta=0.
\end{eqnarray*}
\end{defi}
The map $\Delta$ is called the cobracket of $L$.
Combining the algebra and coalgebra structures under some compatibility condition, we obtain
the definition of a Lie bialgebra as follows.
\begin{defi}\rm
A Lie bialgebra is a triple $(L,\phi,\Delta)$ where $(L,\phi)$ is a Lie algebra, $(L,\Delta)$ is a Lie coalgebra and
\begin{equation}\label{bi1}
\Delta \phi(x\otimes y)= x\cdot\Delta y-y\cdot\Delta x,  \ \ \forall\  x, y\in L.
\end{equation}
\end{defi}
Note that the compatibility condition (\ref{bi1}) is equivalent to say that $\Delta$ is a derivation.
\begin{defi}\rm
A coboundary Lie bialgebra is a $(L,\phi,\Delta,r)$ where $(L,\phi,\Delta)$ is a Lie bialgebra and $\r\in {\rm Im}(1\otimes 1-\tau)$
such that $\Delta$ is a coboundary of $r$, i.e., for arbitrary $x\in L$,
\begin{equation*}
\Delta(x)=x\cdot r.
\end{equation*}
\end{defi}

Denote $U(L)$ the universal enveloping algebra of $L$ and $1$ the identity element
of $U(L)$. For any $r = \sum_i a_i\otimes b_i\in L\otimes L$, define $r^{ij}$ and $c(r)$ to be
the elements of $U(L)\otimes U(L)\otimes U(L)$ by
\begin{equation*}
r^{12}=\sum_i a_i\otimes b_i \otimes 1;\ \
r^{13}=\sum_i a_i\otimes 1\otimes b_i;\ \
r^{23}=\sum_i 1\otimes a_i\otimes b_i;
\end{equation*}
and
$$c(r)= [r^{12}, r^{13}]+[r^{12}, r^{23}]+[r^{13}, r^{23}].$$
The following result is due to Drinfeld \cite{D2}.
\begin{theo}
Let $(L,\phi)$ be a Lie algebra, then $\Delta=\Delta_r$ $($for some $r\in{\rm Im}(1-\tau){\ssc\,})$
endows $(L,\phi,\Delta)$ with a Lie bialgebra structure if and only if r satisfies the following
Modified Yang-Baxter Equation {\rm(}MYBE{\rm)} \rm {:}
\begin{equation*}
x\cdot c(r)=0, \ \ \forall\  x\in L.
\end{equation*}
\end{theo}
Furthermore, one has the following definition.
\begin{defi}\rm
A coboundary Lie bialgebra $(L,\phi,\Delta,r)$ is called triangular if $r$
satisfies the following Classical Yang-Baxter Equation (CYBE) \rm {:}
\begin{equation*}
c(r) = 0.
\end{equation*}
\end{defi}
\section{Lie bialgebras of generalized loop Virasoro algebras}

Let $\Gamma$ be an additive subgroup of $\mathbb{F}$ and $\mathbb{F}[\Gamma]$ be the group algebra of $\Gamma$
with basis $\{e^{\alpha}\,|\,\alpha\in \Gamma\}$.  Denote $\mathscr{W}(\Gamma)=\mathbb{F}[\Gamma]\partial$ (where
$\partial(e^{\alpha})=\alpha e^{\alpha}$), which is usually referred to as the
{\it $($generalized$)$ Witt algebra} or {\it $($generalized centerless$)$ Virasoro
algebra} (c.f. \cite{SZ1}).
Since $\mathscr{W}(\Gamma)\cong\mathscr{W}(\Gamma')$ if $\Gamma\cong\Gamma'$ by \cite{SXZ}, and for any
$0\neq\al\in \Gamma$ we have $\Gamma\cong\frac{1}{\al}\Gamma$, replacing $\Gamma$ by
$\frac{1}{\al}\Gamma$ if necessary, we always assume that
$1\in\Gamma$.
We denote ${\L}=\mathscr{W}(\Gamma)\otimes \mathbb{F}[t^{\pm 1}]$, which is referred to as  a
{\it $($generalized$)$ loop Witt algebra} or {\it $($generalized$)$ loop $($centerless$)$
Virasoro algebra}.
We can   write ${\L}$ as
\begin{equation}\label{a1-0}
{\L}=\raisebox{-5pt}{${}^{\, \, \displaystyle\oplus}_{\alpha\in \Gamma}$}{\L}_{\alpha},\ \ \
 {\L}_{\alpha}={\rm span} \{L_{\alpha}t^{i}\mid i\in \mathbb{Z}\},\ \ \ L_{\alpha}t^{i}=e^{\alpha}\partial\otimes t^{i}\mbox{ \ for }\al\in\Gamma,
 \end{equation}
with brackets
\begin{equation}\label{a1}
[L_{\al}t^{i},L_{\beta}t^{j}]=(\beta-\al)L_{\alpha+\beta}t^{i+j},\ \ \forall\, \alpha,\beta\in\Gamma,\ i,j \in \Z.
\end{equation}
We use the notation $L_{\al,i}:=L_{\al}t^{i}$ in the remaining of this section.

Denote by $V$ the tensor product $\L\otimes\L$. Since $\L$ is a $\Gamma$-graded Lie algebra, $V$ admits a natural $\Gamma$-graded $\L$-module structure under the adjoint diagonal action of $\L$. Precisely, $V=\oplus_{\alpha\in \Gamma} V_\al$ where $V_\al$ is spanned by $\{L_{\beta,i}\otimes L_{\r,j}\, |\,  \beta+\r=\al,\  i,j\in\Z\}$.
Denote by ${\rm Der\,}(\L,V)$
and ${\rm Inn\,}(\L,V)$ the vector spaces of all derivations and inner derivations from $\L$ to $V$ respectively.
It is well known that the first cohomology group of $\L$ with coefficients in the module
$V$ is isomorphic to
${\rm Der\,}(\L,V)/{\rm Inn\,}(\L,V)$.
Denote by $({\rm Der\,}(\L,V))_{\al}$ the space of derivations of degree $\al$.
Then we have the following result through some arguments analogous to that in
\cite{SS}.
\begin{lemm}\label{le1}
The space ${\rm Der}(\L,V)$ can be decomposed into
\begin{equation*}
{\rm Der}(\L,V)=\!\raisebox{-5pt}{${}^{\, \, \displaystyle\oplus}_{\alpha\in \Gamma}$}({\rm Der}(\L,V))_\al
\end{equation*}
where $({\rm Der}(\L,V))_\al\subset {\rm Inn}(\L,V)$ for any $\al\neq0$.\QED
\end{lemm}
We can also obtain the following result easily, which is useful in our paper.
\begin{lemm}\label{le2}
Let $D\in {\rm Der}(\L,V)$ be a derivation such that $D(L_{0,0})=0$, then $D\in ({\rm Der}(\L,V))_0$.
\end{lemm}
\ni{\it Proof.}\ \  Let $D$ act on the equation $[L_{0,0},L_{\al,i}]=\al L_{\al,i}$, we have
$[L_{0,0}, D(L_{\al,i})]=\al D(L_{\al,i})$. So $D(L_{\al,i})\in V_\al$ and the assertion follows.\QED

Fix $i$ and $j$, let $V_{i,j}={\rm span}\{L_{\al,i}\otimes L_{\beta, j}\, |\, \al,\beta\in\Gamma\}={\rm span}\{(L_{\al}\otimes L_{\beta})x^i y^j \,|\, \al,\beta\in\Gamma\}$.
Then $V\!=\!\oplus_{i,j\in\Z}V_{i,j}$, which is a $\Z^2-$graded space.
Consider the subalgebra $\W\!=\!{\rm span}\{L_{\al,0}\,|\,\al\!\in\!\Gamma\}$ of $\L$.
Then $\W$ is actually isomorphic to the generalized centerless Virasoro algebra
and $V_{i,j}$ is a $\W$-module.
It is easy to show that $V_{i,j}\cong\W\otimes\W$ as $\W$-modules.

\begin{lemm}\label{le3}
For  $D\in{\rm Der}(\L,V)$, there is a vector $v\in V$ such that $(D-D_v)(\W)=0$, where $D_v$ is an inner derivation defined by $D_v(L_{\al,i})=L_{\al,i}\cdot v$ for any $\al\in\Gamma$ and $i\in\Z$.
\end{lemm}
\ni{\it Proof.}\ \  Consider the restriction of $D$ to $\W$, we have $D|_\W\in{\rm Der}(\W,V)$.
We also denote $D|_\W$ by $D$ for convenience.
As a $\W$-module, $V=\oplus_{i,j\in\Z}V_{i,j}$.
Then $D=\sum_{i,j\in\Z} D_{i,j}$ (i.e., \cite{SS}), where $D_{i,j}\in {\rm Der}(\W,V_{i,j})$.
We should remark that this formula holds in the sense that
when $D$ acts on a given vector, there exist only finite nonzero terms on the right side.
Referring to \cite{Ng} or \cite{SS}, one has $H^1(\W,\W\otimes\W)=0$.
So there exists some $v_{i,j}\in V_{i,j}$ such that $D_{i,j}=D_{v_{i,j}}$,
where $D_{v_{i,j}}$ is an inner derivation with respect to $v_{i,j}$.
We will show that the formula $D=\sum_{i,j\in\Z} D_{i,j}$ is just a finite sum in the rest of the proof,
then we can deduce that $D$ is an inner derivation.

If $\Gamma$ is a finitely generated free abelian group (e.g., $\Z$, $\Z^2$),  $\W$ is finitely generated as a Lie algebra.
Then we can deduce that $D(\W)$ is contained in a finite direct sum of $V_{i,j}$.
Thus $D=\sum_{i,j\in\Z} D_{i,j}$ is a finite sum and $D=D_v$ where $v=\sum_{i,j\in\Z} v_{i,j}$.

For an arbitrary group $\Gamma$, we have to use completely different arguments.
Consider the finite
set $X=\{(i,j)\,|\,(L_{0,0}\cdot v_{i,j},L_{1,0}\cdot v_{i,j},L_{2,0}\cdot v_{i,j})\neq(0,0,0)\}$.
The complement of $X$ in $\Z\times\Z$ is $\overline{X}=\{(i,j)\,|\,L_{0,0}\cdot v_{i,j}=L_{1,0}\cdot v_{i,j}=L_{2,0}\cdot v_{i,j}=0\}$.
Given a pair $(i,j)\in\overline{X}$ and assume $v_{i,j}=\sum_{\al,\beta\in\Gamma} a(\al,\beta) L_\al\otimes L_\beta x^iy^j$, where $a(\al,\beta)\in\F$.
Since
\begin{equation*}
L_{0,0}\cdot v_{i,j}=\sum_{\al,\beta\in\Gamma} (\beta+\al)a(\al,\beta)L_\al\otimes L_\beta x^iy^j=0,
\end{equation*}
we have $(\beta+\al)a(\al,\beta)=0$.
Then we reduce to the case that $v_{i,j}=\sum_{\r\in\Gamma} a(\r) L_\r\otimes L_{-\r} x^iy^j$ for $a(\r)\in\F$.

Since
\begin{equation*}
L_{1,0}\cdot v_{i,j}=\sum_{\r\in\Gamma} L_\r\otimes L_{1-\r} ((-\r-1)a(\r)+(\r-2)a(\r-1))x^iy^j=0,
\end{equation*}
we have
\begin{equation}\label{e2.2}
(\r+1)a(\r)=(\r-2)a(\r-1),
\end{equation}
and
\begin{equation}\label{e2.3}
\r a(\r-1)=(\r-3)a(\r-2).
\end{equation}
Similarly,
\begin{equation*}
L_{2,0}\cdot v_{i,j}=\sum_{\r\in\Gamma} L_\r\otimes L_{2-\r} ((-\r-2)a(\r)+(\r-4)a(\r-2))x^iy^j=0,
\end{equation*}
and we obtain that
\begin{equation}\label{e2.4}
(\r+2)a(\r)=(\r-4)a(\r-2).
\end{equation}

From the linear equations (\ref{e2.2}),\ (\ref{e2.3}) and (\ref{e2.4}), we have $a(\r)=0$.
Then $v_{i,j}=0$ and $\{(i,j)\,|\,v_{i,j}\neq0\}\subset X$. Thus $D=\sum_{i,j\in\Z} D_{i,j}$ is a finite sum and $D=D_v$ where $v=$ $\sum_{i,j\in\Z} v_{i,j}$.
So ${\rm Der}(\W,V)=\oplus_{i,j\in\Z}{\rm Der}(\W,V_{i,j})={\rm Inn}(\W,V)$. This completes the proof.\QED

For the sake of computation, we will adopt the following description of $V$ in the rest of this section. As a vector space,  $V\cong(\W\otimes\W)\otimes\F[x,x^{-1},y,y^{-1}]\cong(\W\otimes\W)[x,x^{-1},y,y^{-1}]$, which is actually the space of formal Laurant polynomials
in two variables with coefficients in $\W\otimes\W$.
Explicitly, this isomorphism maps $L_{\al,i}\otimes L_{\beta,j}$ to $(L_{\al}\otimes L_{\beta})x^iy^j$.
Now the action of $\L$ on $V$ is given by
\begin{equation}\label{e2.1}
L_{\al,i}\cdot((L_{\beta}\otimes L_{\gamma})x^j y^k)=(\beta-\al)(L_{\al+\beta}\otimes L_{\gamma})x^{i+j} y^k +(\gamma-\al)(L_{\beta}\otimes L_{\al+\gamma})x^j y^{i+k}.
\end{equation}

\begin{prop}\label{p1}
Every derivation from $\L$ to $V$ is inner, i.e., $H^1(\L,V)=0$.
\end{prop}
\ni{\it Proof.}\ \   Take a derivation $D$ from $\L$ to $V$. From Lemma \ref{le3},  replacing $D$ by $D-D_v$ for some $v\in V$,
we can assume $D(\W)=0$. Then we have $D\in({\rm Der}(\L,V))_0$ from Lemma \ref{le2}.
Assume that $D(L_{\al,i})=\sum_{\gamma\in\Gamma}(L_\gamma\otimes L_{\al-\gamma})f_{\al,i,\gamma}$,
where $f_{\al,i,\gamma}\in\F[x,x^{-1},y,y^{-1}]$ and $f_{\al,0,\gamma}=0$.

Let $D$ act on (\ref{a1}), we have
\begin{equation}\label{main}\aligned
(\beta-\al)f_{\al+\beta,i+j,\gamma}
=&(\gamma-2\al)x^if_{\beta,j,\gamma-\al}+(\beta-\al-\gamma)y^i f_{\beta,j,\gamma}\\
&-(\gamma-2\beta)x^jf_{\al,i,\gamma-\beta}-(\al-\beta-\gamma)y^j f_{\al,i,\gamma}.
\endaligned
\end{equation}
Setting $\beta=-\al$ in (\ref{main}), we have
\begin{equation}\label{main2}\aligned
(-2\al)f_{0,i+j,\gamma}
=\ &(\gamma-2\al)x^if_{-\al,j,\gamma-\al}+(-2\al-\gamma)y^i f_{-\al,j,\gamma}-\\
&(\gamma+2\al)x^jf_{\al,i,\gamma+\al}-(2\al-\gamma)y^j f_{\al,i,\gamma}.
\endaligned
\end{equation}
Setting $\al=\beta=0$ in (\ref{main}), we have
\begin{equation*}
(x^i-y^i) \gamma f_{0,j,\gamma}=(x^j-y^j) \gamma f_{0,i,\gamma}.
\end{equation*}
If $\gamma\neq0$ we have
\begin{equation}\label{main3}
(x-y)f_{0,i,\gamma}=(x^i-y^i)f_{0,1,\gamma}.
\end{equation}
Setting $\beta=0$ and $i=0$ in (\ref{main}), we have
\begin{equation*}
-\al f_{\al,j,\gamma}=(\gamma-2\al)f_{0,j,\gamma-\al}+(-\al-\gamma)f_{0,j,\gamma}.
\end{equation*}
Substituting this formula into (\ref{main2}), we obtain that
\begin{equation}\label{main4}
\begin{split}
(-2\al^2)f_{0,i+j,\gamma}
=\ &(\gamma-2\al)x^i((\gamma+\al)f_{0,j,r}-(\gamma-2\al)f_{0,j,\gamma-\al})\\
&+(-2\al-\gamma)y^i((\gamma+2\al)f_{0,j,\gamma+\al}-(\gamma-\al)f_{0,j,\gamma})\\
&+(\gamma+2\al)x^j((\gamma-\al)f_{0,i,\gamma}-(\gamma+2\al)f_{0,i,\gamma+\al})\\
&+(2\al-\gamma)y^j((\gamma-2\al)f_{0,i,\gamma-\al}-(\gamma+\al)f_{0,i,\gamma}).
\end{split}
\end{equation}
If $\al\neq\gamma$, $\al\neq-\gamma$ and $\gamma\neq0$, multiplying $x-y$ on both sides of the above equation we have
\begin{equation}
2(\gamma^2-\al^2)f_{0,1,\gamma}=(\gamma+2\al)^2f_{0,1,\al+\gamma}+(\gamma-2\al)^2f_{0,1,\gamma-\al}.
\end{equation}
Since $D(L_{0,i})=\sum_{\gamma\in\Gamma}(L_\gamma\otimes L_{\al-\gamma})f_{0,i,\gamma}$, there are only finite nonzero terms on the right side.
We can choose some $\al$ such that $\gamma^2-\al^2\neq0$, $f_{0,1,\al+\gamma}=0$ and $f_{0,1,\gamma-\al}=0$, then we deduce that $f_{0,1,\gamma}=0$ if $\gamma\neq0$. Consequently, $f_{0,i,\gamma}=0$ for any $\gamma\neq0$ from (\ref{main3}).

Now we reduce to the case that $D(L_{0,i})\in (L_0\otimes L_0)[x,x^{-1},y,y^{-1}]$.
Assume that $D(L_{0,i})=L_0\otimes L_0 f_i.$
Then $D(L_{\al,i})=L_{\al}\otimes L_0 f_i+ L_0\otimes L_\al f_i$ for any $\al\neq0$ from (\ref{main4}).
Setting $\gamma=\al$ in (\ref{main2}), we get
\begin{equation}
x^if_j+y^jf_i=0,\ \ \forall\ i,j\in\Z.
\end{equation}
Hence we have $f_i=0$ for any $i\in\Z$. Thus $D=0$. \QED
\begin{prop}\label{p2}
Suppose $r\in V$ satisfying $a\cdot r\in {\rm Im}(1\otimes 1-\tau)$ for all $a\in \L$. Then
$r\in {\rm Im}(1\otimes 1-\tau)$.
\end{prop}
\ni{\it Proof.}\ \   Assume that $r=\sum_{\al\in\Gamma} r_\al$. Since $1\otimes 1-\tau$ is a homogenous operator of degree $0$, then $a\cdot r\in {\rm Im}(1\otimes 1-\tau)$ implies $a\cdot r_\al\in {\rm Im}(1\otimes 1-\tau)$. Since $L_{0,0}\cdot r_\al=\al r_\al$, we can obtain that $r_\al\in {\rm Im}(1\otimes 1-\tau)$ for any $\alpha\neq0$.
So $r=r_0=\sum_{\al\in\Gamma} L_\al\otimes L_{-\al} f_\al(x,y)$, where $f_\al(x,y)\in\F[x,x^{-1},y,y^{-1}]$.
Since $L_{0,i}\cdot r=\sum_{\al\in\Gamma} \al L_\al\otimes L_{-\al}(x^i-y^i)f_\al(x,y)\in {\rm Im}(1\otimes 1-\tau)$.
Note that ${\rm Im}(1\otimes 1-\tau)={\rm Ker}(1\otimes 1+\tau)$ and $\tau(L_\al\otimes L_\beta f(x,y))=L_\beta\otimes L_\al f(y,x)$.
Then $\al(x^i-y^i)f_\al(x,y)-\al(y^i-x^i)f_{-\al}(y,x)=0$.
So $f_\al(x,y)+f_{-\al}(y,x)=0$ for any $\alpha\neq0$.
Then $L_\al\otimes L_{-\al} f_\al(x,y)+L_{-\al}\otimes L_\al f_{-\al}(x,y) \in {\rm Im}(1\otimes 1-\tau)$ if $\alpha\neq0$.
It reduces to the case that $r=L_0\otimes L_0 f_0(x,y)$.
Then $L_{\al,0}\cdot r=-\al L_\al\otimes L_0 f_0(x,y)-\al L_0\otimes L_\al f_0(x,y)\in {\rm Im}(1\otimes 1-\tau))$.
We have $\al(f_0(x,y)+f_0(y,x))=0$ for any $\al\in\Gamma$.
Then we obtain that $r\in{\rm Im}(1\otimes 1-\tau)$. This completes the proof.\QED

\begin{prop}\label{p3}
Suppose that a vector $c\in \L\otimes\L\otimes\L$ satisfies $a\cdot c=0$ for all $a\in\L$, then $c=0$.
\end{prop}
\ni{\it Proof.}\ \   Assume $c=\sum_{\al,\beta,\r\in\Gamma} L_\al\otimes L_\beta\otimes L_\r f_{\al,\beta,\r}(x,y,z)$, where $f_{\al,\beta,\r}(x,y,z)\in\F[x^{\pm1},y^{\pm1},z^{\pm1}]$.
Choose an order on $\Gamma$ compatible with the group structure on $\Gamma$.
Then we have an induced lexicographic order on $\Gamma\times\Gamma\times\Gamma$.
Let $(\al_0,\beta_0,\r_0)={\rm max}\{(\al,\beta,\r)\,|\,L_\al\otimes L_\beta\otimes L_\r f_{\al,\beta,\r}\neq0\}$.
Take a $\delta\in\Gamma\setminus\{0\}$ such that $\delta\neq\al_0$. Then the maximal term of $L_{\delta,0}\cdot c$ is $(\al_0-\delta)L_{\al_0+\delta}\otimes L_{\beta_0}\otimes L_{\r_0}f_{\al_0,\beta_0,\r_0}$. Then $f_{\al_0,\beta_0,\r_0}=0$. This completes the proof.\QED

\begin{theo}\label{th1}
Every Lie bialgebra structure on $\L$ is a coboundary triangular Lie bialgebra.
\end{theo}
\ni{\it Proof.}\ \   Suppose $(L,\phi,\Delta)$ is a Lie bialgebra structure on $\L$.
Then $\Delta=\Delta_r$ for some $r\in L\otimes L$.
Since ${\rm Im}\Delta \subset {\rm Im}(1\otimes 1-\tau)$, we get $a\cdot r\in {\rm Im}(1\otimes 1-\tau)$ for all $a\in \L$.
From Proposition \ref{p2}, we have $r\in {\rm Im}(1\otimes 1-\tau)$,
and $c(r)=0$ from Proposition \ref{p3}. Thus $(L,\phi,\Delta)$ is a coboundary triangular Lie bialgebra.\QED


\section{Lie bialgebras of generalized map Virasoro algebras}
In this section, we will generalize the results in Section 3 to the case of map Virasoro algebras.
Let $A$ be a unital commutative associative algebra such that $A\otimes A$ is an integral domain.
We denote ${\mathscr{W}_{A}(\Gamma)}=\mathscr{W}(\Gamma)\otimes A$, which is referred to as  a
{\it $($generalized$)$ map Witt algebra} or {\it $($generalized$)$ map $($centerless$)$
Virasoro algebra}. In case $A={\mathbb{F}}[t,t^{-1}]$, we get the (generalized) loop Virasoro algebras.
We can  write ${{\mathscr{W}_{A}(\Gamma)}}$ as
\begin{equation}\label{me1}
{{\mathscr{W}_{A}(\Gamma)}}\!=\!\raisebox{-5pt}{${}^{\, \, \displaystyle\oplus}_{\alpha\in \Gamma}$}{{\mathscr{W}_{A}(\Gamma)}}_{\alpha},\ \
\mbox{where \ }
{{\mathscr{W}_{A}(\Gamma)}}_{\alpha}\!\!={\rm span} \{L_{\alpha}x \mid \alpha\in\Gamma, x\in A\},
\end{equation}
with relations
\begin{equation}\label{me2}
[L_{\alpha}x,L_{\beta}y]=(\beta-\alpha)L_{\alpha+\beta}xy,\ \ \forall\, \alpha,\beta\in\Gamma,\ x,y \in A.
\end{equation}
Denote the twist map on $A\otimes A$ by $\iota(x\otimes y)=y\otimes x$.

For an element $x\in A$, we denote the left action of $x$ on $A\otimes A$ by $L_x$ and the right action of $x$ on $A\otimes A$ by $R_x$. Explicitly,
$L_x(y\otimes z)=(xy)\otimes z$ and $R_x(y\otimes z)=y\otimes zx$. 
Then we get two different $A-$module structures on $A\otimes A$.

Let $V$ be the ${\mathscr{W}_{A}(\Gamma)}-$module ${\mathscr{W}_{A}(\Gamma)}\otimes {\mathscr{W}_{A}(\Gamma)}$.
As a vector space, $V\cong(\W\otimes\W)\otimes(A\otimes A)$.
Explicitly, this isomorphism maps $L_\al x\otimes L_\beta y$ to $(L_{\al}\otimes L_{\beta})\otimes(x\otimes y)$,
then the action of $\L$ on $V$ is given by
\begin{equation}
L_\al x\cdot((L_{\beta}\otimes L_{\gamma})(y\otimes z))=(\beta-\al)(L_{\al+\beta}\otimes L_{\gamma})(xy\otimes z)+(\gamma-\al)(L_{\beta}\otimes L_{\al+\gamma}) (y\otimes xz).
\end{equation}
Given an element $f\in(A\otimes A)$,
\begin{equation}
L_\al x\cdot((L_{\beta}\otimes L_{\gamma})f)=(\beta-\al)(L_{\al+\beta}\otimes L_{\gamma})(L_xf)+(\gamma-\al)(L_{\beta}\otimes L_{\al+\gamma})(R_xf).
\end{equation}
The twist map of $V$ is given by
\begin{equation}
\tau((L_{\beta}\otimes L_{\gamma})f)=(L_{\gamma}\otimes L_{\beta})\iota(f) .
\end{equation}
Similar to the loop case we have the following result.
\begin{lemm}\begin{enumerate}\item[\rm(1)]
The space ${\rm Der}(\M,V)$ can be decomposed into
\begin{equation*}
{\rm Der}(\M,V)=\!\raisebox{-5pt}{${}^{\, \, \displaystyle\oplus}_{\alpha\in \Gamma}$}({\rm Der}(\M,V))_\al
\end{equation*}
where $({\rm Der}(\M,V))_\al\subset Inn(\M,V)$ for any $\al\neq0$.
\item[\rm(2)]
For any $D\in {\rm Der}(\M,V)$ such that $D(L_0)=0$, we have $D\in ({\rm Der}(\M,V))_0$.
\end{enumerate}\end{lemm}
Suppose $\{x_i\,|\,i\in I\}$ is a basis of $A$.
Then $\{x_i\otimes x_j\,|\,i, j\in I\}$ is a basis of $A\otimes A$.
By a similar argument, we have the following lemma.
\begin{lemm}
For $D\in {\rm Der}(\M,V)$, there exists a vector $v\in V$ such $(D-D_v)(\W)=0$, where $D_v$ is defined by $D_v(L_\al x)=(L_\al x)\cdot v$.\QED
\end{lemm}

Now we can prove the following result.
\begin{prop}\label{p4}
Every derivation from $\M$ to $V$ is inner, i.e., $H^1(\M,V)=0$.
\end{prop}
\ni{\it Proof.}\ \    Similarly, we take a derivation $D\in ({\rm Der}(\M,V))_0$ such that $D(\W)=0$.
Assume that $D(L_\al x)=\sum_{\r\in\Gamma} L_\r\otimes L_{\al-\r}f_{\al,x,\r}$ where  $f_{\al,x,\r}\in A\otimes A$.
Obviously, $f_{\al,1,\r}=0$.
Applying $D$ to (\ref{me2}), one has
\begin{equation}\label{me}
\aligned
(\beta-\al)f_{\al+\beta,xy,\gamma}
=\ &(\gamma-2\al)L_xf_{\beta,y,\gamma-\al}+(\beta-\al-\gamma)R_x f_{\beta,y,\gamma}\\
&-(\gamma-2\beta)L_yf_{\al,x,\gamma-\beta}-(\al-\beta-\gamma)R_y f_{\al,x,\gamma}.
\endaligned
\end{equation}
Letting $\al=\b=0$ in (\ref{me}), we have
\begin{equation}
\r(L_x-R_x)f_{0,y,\r}=\r(L_y-R_y) f_{0,x,\r}.
\end{equation}
Letting $\al+\b=0$ in (\ref{me}), we have
\begin{equation}\label{me3}
\aligned
(-2\al)f_{0,xy,\gamma}
=\ &(\gamma-2\al)L_xf_{-\al,y,\gamma-\al}+(-2\al-\gamma)R_x f_{-\al,y,\gamma}\\
&-(\gamma+2\al)L_yf_{\al,x,\gamma+\al}-(2\al-\gamma)R_y f_{\al,x,\gamma}.
\endaligned
\end{equation}
Letting $x=1$ and $\b=0$ in (\ref{me}), we have
\begin{equation}\label{me4}
\al f_{\al,y,\r}=(\r+\al)f_{0,y,\r}-(\r-2\al)f_{0,y,\r-\al}.
\end{equation}
From (\ref{me3}) and (\ref{me4}), we obtain that
\begin{equation}
\aligned
(-2\al^2)f_{0,xy,\gamma}
=\ &(\gamma-2\al)L_x((\gamma+\al)f_{0,y,r}-(\gamma-2\al)f_{0,y,\gamma-\al})\\
&+(-2\al-\gamma)R_x((\gamma+2\al)f_{0,y,\gamma+\al}-(\gamma-\al)f_{0,y,\gamma})\\
&+(\gamma+2\al)L_y((\gamma-\al)f_{0,x,\gamma}-(\gamma+2\al)f_{0,x,\gamma+\al})\\
&+(2\al-\gamma)R_y((\gamma-2\al)f_{0,x,\gamma-\al}-(\gamma+\al)f_{0,x,\gamma}).
\endaligned
\end{equation}
If $\al\neq\gamma$, $\al\neq-\gamma$, $\gamma\neq0$ and $xy\otimes1-1\otimes xy\neq0$, applying $L_z-R_z$ to the above equation we have
\begin{equation}
2(\gamma^2-\al^2)f_{0,z,\gamma}=(\gamma+2\al)^2f_{0,z,\al+\gamma}+(\gamma-2\al)^2f_{0,z,\gamma-\al}.
\end{equation}
Then $f_{0,z,\gamma}=0$ for any $\r\neq0$.

Now we assume that $D(L_0z)=L_0\otimes L_0f_z$ for some $f_z\in A\otimes A$.
Then $D(L_\al z)=L_\al\otimes L_0f_z+L_0\otimes L_\al f_z$ for any $\al\neq 0$.
Let $D$ act on the equation $[L_{-\al},L_{\al}x]=2 \al L_0 x$, we have $f_x=0$.
Thus $D=0$. This completes the proof.\QED

\begin{prop}\label{p5}
Suppose $r\in V$ such that $a\cdot r\in {\rm Im}(1\otimes 1-\tau)$ for all $a\in \M$. Then
$r\in {\rm Im}(1\otimes 1-\tau)$.
\end{prop}
\ni{\it Proof.}\ \   Similar to the proof of Proposition \ref{p2}, we can take $r=r_0=\sum_{\al\in\Gamma} L_\al\otimes L_{-\al} f_\al$, where $f_\al\in A\otimes A$.
Let $Y=\{\al \,|\,f_\al\neq0\}$ and take some $\b$ not belonging to $Y$,
then $L_\b\cdot r=$ $\sum_{\al\in Y} ((\al-\b)L_{\al+\b} \otimes L_{-\al} f_\al + (-\al-\b) L_\al\otimes L_{-\al+\b} f_\al)$.
Thus $f_\al+\iota(f_{-\al})=0$ and $r\in{\rm Im}(1\otimes 1-\tau)$.
This completes the proof.\QED

Similar to the loop case, we can obtain the following proposition.
\begin{prop}\label{p6}
Suppose $c\in \M\otimes\M\otimes\M$ satisfies $a\cdot c=0$ for all $a\in\M$,
then $c=0$.\QED
\end{prop}

Finally we get the main result in this section.

\begin{theo}\label{th2}
Every Lie bialgebra structure on $\M$ is ar coboundary triangula Lie bialgebra.
\end{theo}
\ni{\it Proof.}\ \  It follows from Propositions \ref{p4}, \ref{p5} and \ref{p6} immediately.\QED

\end{document}